\documentclass[12pt]{article}

\usepackage[a4paper]{geometry}

\usepackage{latexsym}
\usepackage{amssymb}
\usepackage{euscript}
\usepackage{mathrsfs}
\usepackage{euscript}
\usepackage{xspace}
\usepackage{url}
\usepackage{tabularx}
\usepackage{booktabs}
\usepackage{graphicx}
\usepackage{xcolor}
\usepackage{enumerate}
\usepackage{soul}\sethlcolor{lightgray}

\usepackage{stackengine,scalerel}
\usepackage{color}

\usepackage[sort,longnamesfirst]{natbib}
\setcitestyle{authoryear,round,citesep={,},aysep={,},yysep={,},notesep={, }}

    \usepackage{tikz}
    \usetikzlibrary{patterns}

\usepackage[pagewise]{lineno}
\input{ntheoremdefs}
\AtBeginDocument{}

\newcommand{\NN}{\Nat}

\newcommand{\Nat}{\mathbb{N}}

\newcommand{\wrt}{w.r.t.\ }
\newcommand{\resp}{resp.\ }
\newcommand{\ie}{i.e.,\xspace}
\newcommand{\eg}{e.g.,\xspace}
\newcommand{\wl}{w.l.o.g.\xspace}

\newcommand{\D}{\ensuremath{(\D1)}}

\makeatletter
\let\@fnsymbol\@alph
\makeatother

\usepackage{dsfont}

\newcommand{\epr}{\hfill $\Box$\mbox{}\\ }
\newcommand{\bpr}{\noindent {\bf Proof.} \hspace{1 em}}

\newcommand{\mm}{\ensuremath{[m_1]\times[m_2]\xspace}}

\newcommand{\domir}{\ensuremath{>^*}}
\newcommand{\strdomir}{\ensuremath{\gg^*}}

\newcommand{\ssk}{\smallskip}
\newcommand{\msk}{\medskip}

\newcommand{\df}{\ensuremath{d_F}}
\newcommand{\dfh}{\ensuremath{d^h_F}}
\newcommand{\dfv}{\ensuremath{d^v_F}}
\newcommand{\dfd}{\ensuremath{d^d_F}}

\newcommand{\un}{\ensuremath{m_1}}
\newcommand{\deu}{\ensuremath{m_2}}

\newcommand{\etri}{\textsc{Electre Tri}\xspace}

\title{About maximal antichains in a product of two chains: A catch-all note%
\,\thanks{Authors are listed alphabetically. They have contributed equally.}
}
\author{Denis Bouyssou\,\thanks{Former Senior Researcher, CNRS, Paris, France,
e-mail: \protect\url{dbouyssou@gmail.com}. }
\and
Thierry Marchant\,\thanks{Ghent University, Department of Data Analysis,
H.\ Dunantlaan, 1, B-9000 Gent, Belgium,
e-mail: \protect\url{thierry.marchant@UGent.be}.}
\and
Marc Pirlot\,\thanks{Universit\'{e} de Mons, rue de Houdain 9, 7000 Mons, Belgium,
e-mail: \protect\url{marc.pirlot@umons.ac.be}.}}

\date{October 2024}

\begin{document}
\maketitle

\section{Introduction}\label{sec:intro}
Our interest for maximal antichains arose from our analysis \citep{BouyssouMP2023TrinBpseudoDisj} of a multicriteria classification method known as the \emph{pseudo-disjunctive \etri-nB} method \citep{FernandezEJOR2017}. We wanted to count the number of maximal antichains in a product of an arbitrary number of finite chains. We did it for products of chains of small cardinality. In the case of a product of two chains of equal sizes, a search in the On-line Encyclopedia of Integer Sequences \citep{OEIS2024} led us to Sequence A171155. Surprisingly, we found no mention of maximal antichains in the various mathematical structures mentioned on the site as counted by this sequence. This motivated us to establish the links between maximal antichains and some of these structures.

\ssk This is what we do in Section~\ref{sec:correspondences}, which establishes one-to-one correspondences (bijections) between maximal antichains in products of two finite linear orders and other mathematical objects.

\ssk Section~\ref{sec:counting} gathers what is known regarding the number of maximal antichains in products of two finite linear orders and establishes new results.

\ssk Some proofs and numerical results are deferred to an Appendix.

\ssk The rest of the present section introduces definitions and some preliminary results. We refer to \citet{CaspardBookEng2012} for notions related to order relations.

\subsection{Dominance orders}
Let $[m]$ denote the integer set $\{1, \ldots, m\}$. We consider the sets $[m_1]$ and $[m_2]$ linearly ordered by the natural order $\geq$ on the integers $\mathbb{N}$. We use $>, < , \leq $ as expected.

\ssk Two natural partial orderings can be defined on the product set $[m_1] \times [m_2]$.
\begin{description}
  \item[The dominance order $\domir$] is defined as $(x,y) \domir (z,w)$ iff $x \geq z$,  $y \geq w$, and at least one inequality is strict; it is an irreflexive, asymmetric and transitive relation.
  \item[The strong dominance order $\strdomir$] is defined as $(x,y) \strdomir (z,w)$ iff $x > z$ and  $y > w$; it is also an irreflexive, asymmetric and transitive relation.
\end{description}

\subsubsection{Chains and antichains}
A \emph{chain} in $\mm$ is a set of distinct elements $\{(x_i, y_i), i \in I\}$, where $ I = [k]$ for some integer $k$, such that, $\forall i \neq j \in I$, $x_i > x_j$ entails $y_i \geq y_j$. It is thus a chain  \wrt the dominance order $\domir$. In other words, $\forall i,j \in I$, we may not have $x_i > x_j$ and $y_i < y_j$.

\smallskip  Note that, \wl, we may number the elements of a chain in a way that respects the lexicographic order on the elements, \ie $x_i > x_{i-1}$ or $[x_i = x_{i-1} \ \textrm{and } \ y_i > y_{i-1}]$, for all $i$.

\smallskip A \emph{strict chain} $\{(x_i, y_i), i \in I\}$ is a chain \wrt the strong dominance order $\strdomir$, \ie $\forall i \neq j \in I$, we have $[x_i > x_j$ and $y_i > y_j]$ or $[x_i < x_j$ and $y_i < y_j]$. This implies that $[x_i \neq x_j$ and $y_i \neq y_j]$ for all $i \neq j$.

\ssk We may number the elements of a strict chain in such a way that $x_i > x_j$, for all $i >j$.

\ssk An \emph{antichain} in $\mm$ is a set of elements $\{(z_i, w_i), i \in I \}$ such that, for all $i\neq j$, we have either $[z_i > z_j$ and $w_i < w_j]$ or $[z_i < z_j$ and $w_i > w_j]$. In other words, all pairs of elements in an antichain are incomparable \wrt the dominance order \domir. In an antichain $\{(z_i, w_i), i \in I\}$, for all $i \neq j$, we have $z_i \neq z_j$ and $w_i \neq w_j$.

A notion of \emph{weak antichain} can be dually defined. The set $\{(z_i, w_i), i \in I\}$ is a weak antichain if its elements are incomparable \wrt  the strong dominance order $\strdomir$. Explicitly: for all $i \neq j$, if $z_i > z_j$ then $w_i \leq w_j$.

\ssk A \emph{maximal antichain} (\resp \emph{maximal strict chain}, \emph{maximal weak antichain}, \emph{maximal chain} is an antichain (\resp strict chain, weak antichain, chain) that is not properly included in another antichain (\resp strict chain, weak antichain, chain).

\section{Correspondences}\label{sec:correspondences}

In this section, we establish links between (maximal) chains and antichains in products of two finite chains and other mathematical structures mentioned in relation with sequence A171155 from \citet{OEIS2024}. We slightly more generally deal with products of two finite chains possibly of unequal lengths $\mm$.

\subsection{Antichains and strict chains}
We have the following correspondence (in products of two chains).

\begin{proposition}\label{prop:bijectStrChainAntichain}
There is a bijection between strict chains and antichains in \mm.
\end{proposition}

\bpr
If $\{(x_i, y_i), i \in I\}$ is a strict chain in \mm, $\{(x_i, m_2 + 1 -y_i), i \in I\}$ is an antichain. Conversely, if $\{(z_i, w_i), i \in I\}$ is an antichain in \mm, $\{(z_i, m_2 + 1 -w_i), i \in I\}$ is a strict chain.
\epr

\ssk This bijection also defines a bijection between maximal strict chains and maximal antichains.

\ssk In a similar way, there is a bijection between chains (\resp maximal chains) and weak antichains (\resp maximal weak antichains) in \mm.

\subsection{Noses of a step matrix and augmentation matrix}
Any antichain in $\mm$ can be associated with two $m_1 \times m_2$ step matrices. Consider for example the antichain $\{(2,4), (4,2)\}$ in $[5]\times[6]$. We represent these two pairs by a boldface ``0'' in the matrix in Table \ref{ta:matrixAntichDominated}. We also represent by a ``0'' in the matrix all pairs that are $\domir$-dominated by an element in the antichain. The remaining pairs can be assigned the value ``1'' (not represented in Table \ref{ta:matrixAntichDominated}). This matrix is a step matrix (whether looking at the 0's or to the 1's). The elements in the antichain, \ie the boldface ``0'' pairs, have been called the \emph{noses} of the ``0'' step matrix \citep[][p.~77]{pirlotvincke1997semiorders}. They are the non $\domir$-dominated pairs in the set of ``0''-valued pairs in the matrix.

\begin{table}[h!!!]
\begin{center}
\begin{tabular}{|c||c|c|c|c|c|c|}
  \hline
    & 1 & 2 & 3 & 4 & 5 & 6 \\
    \hline \hline
  1 & 0 & 0 & 0 & 0 &  &  \\ \hline
  2 & 0 & 0 & 0 & $\mathbf{0}$  &  &  \\ \hline
  3 & 0 & 0 &  &  &  &  \\ \hline
  4 & 0 & $\mathbf{0}$ &  &  &  &  \\ \hline
  5 &  &  &  &  &  &  \\
  \hline
\end{tabular}
\end{center}
\caption{The antichain $\{(2,4), (4,2)\}$ in $[5]\times[6]$ and the elements in $\mm$ that are $\domir$-dominated by an element in the antichain are represented by a ``0'' value.}
\label{ta:matrixAntichDominated}
\end{table}

\ssk The properties of such a matrix representation are general. The ``0'' corresponding to the elements of an antichain and to the elements that are $\domir$-dominated by one of them form a step matrix. More precisely, it is a North-West step matrix (NW step matrix, in the sequel) The elements of the antichain correspond to the noses of this step matrix, since they are maximal \wrt $\domir$ in the set of elements that do not dominate an element in the antichain.

\ssk A second step matrix can be associated with an antichain by writing a ``0'' in all positions corresponding either to an element of the antichain or an element $\domir$-dominating an element of the antichain.  For the example of the antichain $\{(2,4), (4,2\})$ in $[5]\times[6]$, we get the matrix in Table \ref{ta:matrixDominateAntich}.

\begin{table}[h!!!]
\begin{center}
\begin{tabular}{|c||c|c|c|c|c|c|}
  \hline
    & 1 & 2 & 3 & 4 & 5 & 6 \\
    \hline \hline
  1 &  &  &  &  &  &  \\ \hline
  2 &  &  &  & $\mathbf{0}$  & 0 & 0 \\ \hline
  3 &  &  &  & 0 & 0 & 0 \\ \hline
  4 &  & $\mathbf{0}$ & 0 & 0 & 0 & 0 \\ \hline
  5 &  & 0 & 0 & 0 & 0 & 0 \\
  \hline
\end{tabular}
\end{center}
\caption{The antichain $\{(2,4), (4,2)\}$ in $[5]\times[6]$ and the elements in $\mm$ that $\domir$-dominate an element in the antichain are represented by a ``0'' value.}
\label{ta:matrixDominateAntich}
\end{table}

\ssk This matrix is also a step matrix. It is a South-East step matrix (SE step matrix in the sequel).  The elements in the antichain correspond to the noses of this matrix (the minimal elements \wrt $\domir$ in the set of elements dominating an element in the antichain).

\subsubsection{Augmentation matrix of an antichain}\label{sssec:augmentMatrixAntich} Putting together the two step matrices associated with an antichain, \ie writing a 0 in a new matrix whenever there is a 0 in one of the two step matrices, we get a representation of the elements that cannot be added to the antichain because they are either dominated or dominate an element in the antichain. 

\ssk We call \emph{augmentation matrix} associated with an antichain in $\mm$, the binary $m_1 \times m_2$  matrix obtained by assigning the value 0 to the elements in the antichain and those that $\domir$-dominate or are $\domir$-dominated by an element in the antichain; the value 1 is assigned to the remaining positions in the matrix\footnote{We always assume that the matrix rows (\resp columns) are labelled by the elements of $[m_1]$ (\resp $[m_2]$) in increasing natural order.}. The 1's in the matrix, if any, correspond to pairs in $\mm$ that can be added to the antichain (one at a time), making it a larger antichain.

\ssk In the example of the antichain $\{(2,4), (4,2)\}$ in $[5]\times[6]$ above, this yields the augmentation matrix in Table \ref{ta:matrixAugmentAntich}. This matrix indicates that [$(1,5)$ or $(1,6)$] and $(3,3)$ and $(5,1)$ can be added to the antichain without losing the antichain property.

\begin{table}[h!!!]
\begin{center}
\begin{tabular}{|c||c|c|c|c|c|c|}
  \hline
    & 1 & 2 & 3 & 4 & 5 & 6 \\
    \hline \hline
  1 & 0 & 0 & 0 & 0 &  &  \\ \hline
  2 & 0 & 0 & 0 & $\mathbf{0}$ & 0 & 0 \\ \hline
  3 & 0 & 0 &  & 0 & 0 & 0 \\ \hline
  4 & 0 & $\mathbf{0}$ & 0 & 0 & 0 & 0 \\ \hline
  5 &  & 0 & 0 & 0 & 0 & 0 \\
  \hline
\end{tabular}
\end{center}
\caption{Augmentation matrix of the antichain $\{(2,4), (4,2)\}$ in $[5]\times[6]$. 
The non-zero positions correspond to elements that can be added to the antichain while keeping the antichain property.}
\label{ta:matrixAugmentAntich}
\end{table}


\ssk Let us recall the following definition \citep[see \eg][]{FulkersonGross1965IncidenceMA}. A 0-1 matrix has the \emph{consecutive ones property} in rows (\resp columns) if, in each row (\resp column), there is no 0 between two 1's.
\begin{proposition}\label{prop:Consecutive1Prop}
The augmentation matrix of an antichain has the consecutive ones property in rows and columns.
\end{proposition}

\bpr
The fact that ones in each row (\resp column) of the matrix appear all together (\ie without 0 between two ones) results from the fact that a one appears in the augmentation matrix when there is a one at the same position in the two step matrices associated with the antichain. In step matrices, a 0 never appears between two 1's in the same row or the same column.
\epr

\ssk The following characterization of maximal antichains obviously results from the definition of the augmentation matrix.
\begin{proposition}\label{prop:AugmentAntichMax}
An antichain is maximal iff its augmentation matrix is the null matrix.
\end{proposition}

\subsubsection{Step matrices and augmentation matrix for strict chains}\label{sssec:StrChainMatrices} By duality (Proposition \ref{prop:bijectStrChainAntichain}), two step matrices and an augmentation matrix can be associated with any strict chain. For illustration purposes, consider the strict chain $\{(2,3), (4,5)\}$ in $[5]\times[6]$ which is in bijection with the antichain in the examples above.
The first associated step matrix has 0's corresponding to the elements of the strict chain and in the positions associated with elements not better on the first dimension and better on the second than an element in the strict chain. These positions cannot be used to augment the strict chain. The 0's form a North-East (NE) step matrix. The ones form a South-West (SW) step matrix. The second associated step matrix has 0's corresponding to the elements of the strict chain and in the positions associated with elements better on the first dimension and not better on the second than an element in the strict chain. These positions cannot be used to augment the strict chain. The 0's form a SW step matrix. The ones form a NE step matrix. The augmentation matrix associated with the strict chain has a one iff the two step matrices have a one in the same position. Only the elements corresponding with a one in the augmentation matrix can be used to extend the strict chain into a larger strict chain. The augmentation matrix corresponding to the strict chain $(\{2,3), (4,5)\}$ is represented in Table \ref{ta:matrixAugmentStrChain}. Propositions \ref{prop:Consecutive1Prop} and \ref{prop:AugmentAntichMax} transpose immediately to strict chains.

\begin{table}
  \centering
  \begin{tabular}{|c||c|c|c|c|c|c|}
    \hline
     & 1 & 2 & 3 & 4 & 5 & 6 \\ \hline \hline
    1 &  &  & 0 & 0 & 0 & 0 \\ \hline
    2 & 0 & 0 & $\mathbf{0}$ & 0 & 0 & 0 \\ \hline
    3 & 0 & 0 & 0 &  & 0 & 0 \\ \hline
    4 & 0 & 0 & 0 & 0 & $\mathbf{0}$ & 0 \\ \hline
    5 & 0 & 0 & 0 & 0 & 0 &  \\
    \hline
  \end{tabular}
  \caption{Augmentation matrix associated with the strict chain $\{(2,3), (4,5)\}$ in $[5]\times[6]$. Boldface 0's represent the elements in the strict chain. Ones take place in the empty cells. }\label{ta:matrixAugmentStrChain}
\end{table}

\subsection{Alignments of two strings}\label{ssec:alignment}
An \emph{alignment} of a sequence of $m_1$ letters and a sequence of $m_2$ letters is ``a way of pairing up elements of the two strings, optionally skipping some elements but preserving the order'' \citep{Covington2004}.

\ssk Alignments can be represented as follows. Consider for instance a sequence of 4 letters A, B, C, D and a sequence of 3 letters X, Y, Z. An example of an alignment is given in Table \ref{ta:alignmentExpl}.

\begin{table}[h!!!]
  \centering
  \begin{tabular}{cccccc}
  A & $-$ & B & C & $-$ & D \\
  $-$ & X & Y & $-$ & Z & $-$ \\
\end{tabular}
  \caption{An alignment of strings ABCD and XYZ matching only B and Y}\label{ta:alignmentExpl}
\end{table}

\ssk Letters in the same column are matched (or paired up). A letter in the same column as  ``$-$'' is unmatched. The sign ``$-$'' denotes a \emph{skip} \citep{Covington2004}, for obvious reasons. In the example of Table~\ref{ta:alignmentExpl}, only B and Y are matched.

\ssk Another alignment that matches the same letters is represented in Table \ref{ta:alignmentExplVariant}.
\begin{table}[h!!!]
  \centering
  \begin{tabular}{cccccc}
  $-$ & A & B & $-$ & C & D \\
  X & $-$ & Y & Z & $-$ & $-$ \\
\end{tabular}
  \caption{Another alignment of strings ABCD and XYZ also matching only B and Y }\label{ta:alignmentExplVariant}
\end{table}

In the sequel, we shall consider alignments that match the same subsets of letters as equivalent. This corresponds to the ``different middle set'' of alignments in \citet[][Section 7]{Covington2004}.
In other words, we call \emph{assignment} a function $a$ from a subset of $[m_1]$ to a subset of $[m_2]$ respecting the orders on $[m_1]$ and $[m_2]$, \ie if $x, y$ belong to the domain of $a$ and $x < y$, then $a(x) < a(y)$. Note that the latter property implies that the function is injective.

\ssk

\begin{remark}\label{rem:conventionAlign} With such a definition of an assignment, note that the only thing that matters for describing an alignment is the set of matched pairs of letters. The representation in Table \ref{ta:alignmentExpl} is conventional: one starts by putting a ``$-$'' sign in the column of unmatched letters in the first string preceding a matched pair, then we do the same for unmatched letters in the second string. Adopting this convention of representation, we may thus unequivocally describe an alignment by the assignment function $a$, \ie by listing the pairs it matches. For instance, the alignment in Table~\ref{ta:alignmentExpl} is $\{\text{(B,Y)}\}$.
\end{remark}

\subsubsection{Alignments and strict chains}\label{sssec:alignAntich}
In the alignment example in Table \ref{ta:alignmentExpl}, the only pair that is matched is (B,Y), the second letter in the string ABCD with the second in the string XYZ. Let us interpret the pair (B,Y) as a strict antichain in the Cartesian product of the set of string letters, ordered alphabetically. We may alternatively code the letters as numbers: $[4]$ coding the letters A, B, C, D; $[3]$ coding the letters X,Y,Z. Since we consider $(2,2) = $ (B,Y) as the only pair in the strict chain, we have that a pair incomparable to $(2,2)$ may not be added to the strict chain; such pairs correspond to illegal matchings given that B and Y are matched. Such pairs are represented by a 0 in the matrix in Table \ref{ta:augmMatrixAlign}. The 0's in this matrix are the union of the 0's in the NE and the SW step matrices associated with the strict chain.

\begin{table}[h!!!]
  \centering
  \begin{tabular}{|c|c||c|c|c|}
    \hline
     &  & X & Y & Z \\ \hline
     &  & 1 & 2 & 3 \\ \hline \hline
    A & 1 &  & 0 & 0 \\ \hline
    B & 2 & 0 & $\mathbf{0}$ & 0 \\ \hline
    C & 3 & 0 & 0 &  \\ \hline
    D & 4 & 0 & 0 &  \\
    \hline
  \end{tabular}
  \caption{Augmentation matrix associated with an alignment (represented by boldface 0). The empty cells correspond to possible positions for augmenting the alignment.}\label{ta:augmMatrixAlign}
\end{table}

\ssk The alignment above is not maximal. It can be augmented by matching A with X and C with Z or D with Z.

\ssk Conversely, we may interpret any strict chain as an alignment. Consider \eg  the strict chain $\{(2,3), (4,5)\}$ in $[5]\times[6]$. Interpreting $[5]$ (\resp $[6]$) as the string ABCDE (\resp UVWXYZ), the strict chain $\{(2,3), (4,5)\}$ can be interpreted as the alignment in Table \ref{ta:alignStrChain}.

\begin{table}[h!!!]
  \centering
  \begin{tabular}{ccccccccc}
    A & $-$ & $-$ & B & C & $-$ & D & E & $-$\\
    $-$ & U & V & W & $-$ & X & Y & $-$ & Z\\
  \end{tabular}
  \caption{Alignment corresponding to the strict chain $\{(2,3), (4,5)\}$}\label{ta:alignStrChain}
\end{table}

\ssk The augmentation matrix in Table~\ref{ta:augmMatrixAlign2} indicates that the alignment $\{$(B,W), (D,Y)$\}$ can be augmented by matching A with U or V, matching C with X, and matching E with Z.

\begin{table}[h!!!]
  \centering
  \begin{tabular}{|c|c||c|c|c|c|c|c|}
    \hline
     & & U & V & W & X & Y & Z  \\ \hline
     & & 1 & 2 & 3 & 4 & 5 & 6 \\ \hline \hline
    A & 1 &  &  & 0 & 0 & 0 & 0\\ \hline
    B & 2 & 0 & 0 & $\mathbf{0}$ & 0 & 0 & 0 \\ \hline
    C & 3 & 0 & 0 & 0 &  & 0 & 0\\ \hline
    D & 4 & 0 & 0 & 0 & 0 & $\mathbf{0}$ & 0  \\ \hline
    E & 5 & 0 & 0 & 0 & 0 & 0 & \\
    \hline
  \end{tabular}
  \caption{Augmentation matrix corresponding with alignment in Table \ref{ta:alignStrChain}}\label{ta:augmMatrixAlign2}
\end{table}

\subsubsection{Alignments without alternate skips and maximal strict chains}\label{sssec:alternatingSkips} \citet{Covington2004} considered different types of restrictions on alignments. Among them, the ``small set'' of alignments consists of alignments without \emph{alternate skips}. Alternate skips are a succession of two skips, denoted by a $-$ sign, one letter skipped in the first string followed by one in the second string. 

\ssk For example, in the alignment in Table \ref{ta:alignStrChain}, A is skipped and U is skipped just after. Similarly, C and X are successively skipped, as well as E and Z. The corresponding alignment thus does not belong to the ``small set''. In terms of strict chains, the existence of alternate skips means that the chain can be augmented by adding the pair of consecutive skipped elements from the two strings; this strict chain is not maximal. Conversely, if a strict chain is not maximal, it can be augmented, which implies the existence of alternate skips. Indeed, by the representation convention of alignments using skips (``$-$''), before the first matching (\resp between any two consecutive matchings, after the last matching), all skipped letters from String~1 are first listed, then all  skipped letters from String~2 are listed. The absence of alternate skips means that only letters of one string are skipped before the first matching (\resp between any two consecutive matchings, after the last matching). This implies that the alignment cannot be augmented. This proves the following proposition.

\ssk
\begin{proposition}\label{prop:augmentAlign}
An alignment can be augmented iff there are alternate skips.
\end{proposition}

\begin{proposition}\label{prop:maxAlignMaxStrChain}
An alignment has no alternate skips iff the corresponding strict chain is maximal.
\end{proposition}

\bpr
By the previous proposition, an alignment is maximal iff it has no alternate skips. Clearly, a maximal alignment corresponds to a maximal strict chain.
\epr

\begin{remark}\label{rem:alignAntich} A direct correspondence between alignments and antichains is obtained by numbering the elements of one of the strings in reverse order. For instance, in the example in Table \ref{ta:alignStrChain}, we may number Z (\resp Y, X, W, V, U) by 1 (\resp 2, 3, 4, 5, 6). The alignment BW and DY thus corresponds to the antichain $\{(2,4), (4,2)\}$.
\end{remark}

\subsection{Words using an alphabet of three letters}\label{ssec:words3}
We consider words composed by using three distinct letters. Let us use the letters $h, v, d$ (for reasons that will become clear later).

\subsubsection{Alignments and words}\label{sssec:alignWords}
An alignment can be bijectively associated with a word based on an alphabet of three letters $h, v, d$. Each column (in the representation illustrated in Table \ref{ta:alignmentExpl}) is coded as a letter: $v$ (\resp $h$) represents a column with a letter in the first (\resp second) row and a skip ``$-$'' in the second (\resp the first) row; letter $d$ represents a column matching two letters. The example in Table \ref{ta:alignmentExpl} is thus represented by the word $vhdvhv$. The one in Table \ref{ta:alignStrChain} is represented by the word $vhhdvhdvh$.

\ssk The alignments corresponding to a word without consecutive $vh$ or $hv$ are called \emph{alignments without alternate skips}; they form the ``small set'' of alignments according to \citet{Covington2004}. They are alignments in which an unmatched letter in the first sequence is not followed by an unmatched letter in the second sequence and vice versa. The above observations are summarized in the following proposition.

\begin{proposition}\label{prop:alignWords}
Alignments in \citeauthor{Covington2004}'s ``different middle set'' are uniquely represented by a word composed of the letters $h,v,d$. The letter $d$ represents a matched pair in the two sequences, while $h$ (\resp $v$) means that a letter from the first (\resp second) string is skipped. In order to guarantee that the representation of an alignment by a word is unique, it is agreed that between two consecutive occurrences of $d$, all the $v$'s come first followed by all the $h$'s; the same holds before the first (\resp after the last) occurrence of a $d$. The number of letters in the word is the sum of the numbers of unmatched letters in the two strings and the number of matched pairs of letters.

\ssk \noindent Alignments in \citeauthor{Covington2004}'s ``small set'' are uniquely represented by a word composed of the letters $h,v,d$ in which the sequences $hv$ and $vh$ never show up.

\end{proposition}

\subsubsection{Strict chains and words via grid lines}\label{sssec:strChainsWords}
\begin{figure}[h!!!]
  \centering
  \includegraphics[width=7cm]{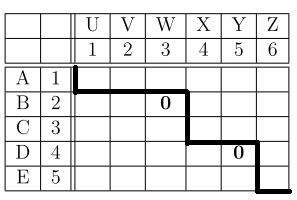}
  \caption{Matrix of strict chain $\{(2,3), (4,5)\}$ represented in Table \ref{ta:augmMatrixAlign2}. The boldface grid line separates the SW region determined by the chain from the other cells. }\label{fig:strChainWord}
\end{figure}

Consider the matrix of a strict chain such as $\{(2,3), (4,5)\}$ in $[5]\times [6]$, represented in Figure \ref{fig:strChainWord}. Looking at the grid containing the elements of the matrix, we aim to describe a grid line starting from the grid's NW corner, ending up in the grid's SE corner, and separating the cells with first coordinate not smaller and second coordinate not larger than those of an element in the chain, \ie the elements in the SW region determined by the chain, from the other cells. We may describe the line as a sequence of vertical downward moves and horizontal rightward moves, which we denote by $v$ and $h$, respectively. The sequence corresponding to the line delimiting the SW region delimited by the strict chain $\{(2,3), (4,5)\}$ in Figure \ref{fig:strChainWord} is $vhhhvvhhvvh$. For turning the ``corner'' determined by an element of the strict chain, the line performs a $h$-move followed by a $v$-move. We substitute all $hv$ sub-sequences with a $d$-move, which can be interpreted as a \emph{diagonal} move crossing the cell containing an element of the strict chain. This yields the following word describing the line represented in Figure \ref{fig:strChainWord}: $vhhdvhdvh$.

\ssk
The following proposition is easy to prove.

\begin{proposition}\label{prop:strChainSWmatrix} Any strict chain in $\mm$ determines a SW step matrix formed of the elements of $\mm$ that have first coordinate not smaller and second coordinate not larger than those of an element in the chain.  Conversely, any SW step matrix corresponds to a strict chain of $\mm$. The elements of the strict chain are the \emph{noses} of the SW step matrix.
\end{proposition}

\ssk The grid line separating the SW step matrix associated with a strict chain is uniquely described by a word composed of the letters $h, v, d$, in which the sequence $hv$ does not appear. The $d$'s correspond to the elements of the strict chain (and the noses of the step matrix). This means that between two consecutive occurrences of $d$'s all $v$'s come first, followed by the $h$'s; the same holds before the first (\resp after the last) occurrence of a $d$.

\ssk
\begin{proposition}\label{prop:maxiStrChainWordSansvh} A strict chain is maximal iff its associated $(h,v,d)$-word does not contain the sequence $vh$.
\end{proposition}

\bpr
If the word associated with a strict chain contains the sequence $vh$, this sequence can be replaced by $d$. The corresponding cell in the matrix is an element that can be added to the chain. Indeed, assume that the sequence $vh$ occurs between two consecutive $d$'s. Let the smallest (\resp largest) of the $d$'s correspond to the chain element of coordinates $(x_1, x_2)$ (\resp $(y_1, y_2)$). The presence of the sequence $vh$ in between implies that $y_1 \geq x_1 + 2$ and $y_2\geq x_2 +2$. Therefore, adding $(x_1 +1, x_2+1)$ to the chain yields a longer chain. A similar reasoning holds in case $vh$ occurs before (\resp after) the first (\resp last) $d$ in the word.

Conversely, if the strict chain is not maximal, an element can be added between two consecutive elements of the chain (or before the first one or after the last one). Using the same reasoning as above, in reverse order, shows that the grid line associated with the chain must contain a $vh$ sequence.
\epr

\subsubsection{Antichains and words via grid lines}\label{sssec:antichWords}

\begin{figure}[h!!!]
  \centering
  \includegraphics[width=7cm]{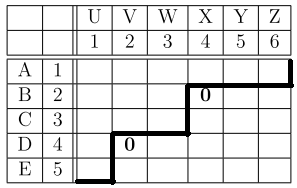}
  \caption{Matrix of antichain $\{(2,4), (4,2)\}$ in $[5]\times[6]$. The boldface grid line separates the SE region determined by the chain from the other cells.}\label{fig:antichWord}
\end{figure}

The grid line associated with an antichain separates the cells either belonging to the antichain or  corresponding to elements $\domir$-dominating an element in the antichain from the rest of the elements of $\mm$. For the example of the antichain $\{(2,4), (4,2)\}$ in $[5]\times [6]$, the grid line is represented in Figure \ref{fig:antichWord}. It separates the SE region determined by the antichain from the complementary NW region. The grid line starting from the NE corner ends up in the SW corner of the grid. It can be described by a $(h, v, d)$-word where $h$ represents a leftward horizontal move, $v$ a downward vertical move, and $d$ an anti-diagonal move. In the case of the example illustrated in Figure \ref{fig:antichWord}, the grid line is represented by the word $vhhdvhdvh$. This is exactly the same word as that representing the dual strict chain $\{(2,3),(4,5)\}$, but the meaning of $h$ and $d$ is different.

\ssk In a similar way as for strict chains, we have the following proposition.

\begin{proposition}\label{prop:antichWords}
An antichain is maximal iff its associated $(h,v,d)$-word does not contain the sequence $vh$.
\end{proposition}


\subsection{Walks on a grid}\label{ssec:antichWalks}
In the previous section, words of two or three letters were associated to strict chains or antichains.
This is done by constructing matrices associated with chains or antichains and considering lines separating some sets of cells in the matrix (grid lines).
In this section, we analyze the relationships between certain walks on a grid (or lattice paths\footnote{Lattice paths have
been studied for a long time \citep[see, e.g.][]{Mohanty80}.
A lattice path
 in the $d$-dimensional integer lattice $\displaystyle \mathbb{Z}^{d}$ with steps in the set $S$,
 is a sequence of vectors  $ v_0 , v_1 ,\ldots , v_k \in \mathbb{Z}^{d}$  such that each consecutive
difference
 $\displaystyle v_{i}-v_{i-1}$ lies in $S$ \citep[][p.21]{Stanley2011}.
Various additional constraints can be imposed on the considered lattice paths.
 For instance, a \emph{Dyck path} is a lattice path starting in $(0,0)$, ending in $(n,n)$,
 composed of $(1,0)$ and $(0,1)$ steps, which never crosses the diagonal $y=x$ but may touch it.
Dyck paths are counted by the Catalan numbers
 \citep[see, e.g.][]{StanleyEnumerativeCombVol2}
 } )
and
strict chains or antichains in a product of chains.
In contrast with the previous section we directly consider grids and certain walks on grids,
without referring to matrices.

\subsubsection{Walks and strict chains}\label{sssec:walksStrChains}
\begin{figure}[h!!!]
  \centering
  \includegraphics[width=7cm]{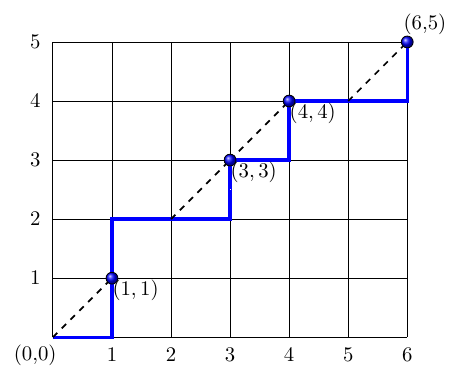}
  \caption{A grid walk from $(0,0)$ to $ (6,5)$. The strict chain $\{(1,1), (3,3), (4,4), (6,5)\}$ is not maximal because it has a $v'h'$ subsequence (namely $(1,1), (1,2), (2,2)$) disjoint from all $h'v'$ subsequences.}\label{fig:gridWalk}
\end{figure}

Consider the integer grid from $(0,0)$ to $(m_1, m_2)$. We are interested in walks starting from $(0,0)$ and ending up in $(m_1, m_2)$ with the property that each step increases one coordinate by one unit while leaving the other coordinate unchanged. A step may thus be denoted by $h'$ (\resp $v'$) if it increases the first (\resp second) coordinate by one unit (we use a different notation in order to avoid confusion with the $h, v$ used before; the relationships between these notations will be analyzed below). For example, Figure~\ref{fig:gridWalk} displays a grid walk from $(0,0)$ to $(6,5)$, which reads $h'v'v'h'h'v'h'v'h'h'v'$. Clearly, the sequence of endpoints of the walk steps forms a chain in $[m_1]\times [m_2]$. In Figure~\ref{fig:gridWalk}, this sequence is $\{(1,0), (1,1), (1,2), (2,2), (3,2), (3,3), (4,3), (4,4), (5,4), (6,4), (6,5)\}$. We have the following result.

\begin{proposition}\label{prop:walkStrictChain}
The endpoints of all subsequences $h'v'$ in a walk from $(0,0)$ to $(m_1, m_2)$ define a strict chain in \mm. Such a strict chain is maximal iff no subsequence $v'h'$ disjoint from the $h'v'$ subsequences appears in the walk .
\end{proposition}
\bpr
Let ${(x_i, y_i), i \in [k]}$, for some integer $k$, be such endpoints numbered in the order in which the corresponding $h'v'$ subsequences occur in the walk.  Each such endpoint $\strdomir$-dominates the previous one and is $\strdomir$-dominated by the next one. If a $v'h'$ subsequence occurs in the walk, then the endpoint of the subsequence $v'h'$ can be added to the strict chain making it a strict chain with more elements. Conversely, if the strict chain ${(x_i, y_i), i \in I}$ is not maximal, an element $(x,y)$ can be added to the chain. The number $x$ is in one of the following three positions: 1) $x < x_1$; 2) $x_i < x < x_{i+1}$ for some $i \in [k-1]$; 3) $x_k < x$. Let us consider the second case (the others are dealt with similarly). Since $x_i < x < x_{i+1}$, we must have $y_i < y < y_{i+1}$, otherwise, $(x,y)$ would be incomparable either to $(x_i, y_i)$ or to $(x_{i+1}, y_{i+1})$. As a consequence, we have $x_{i+1} \geq x_i + 2 $ and $y_{i+1} \geq y_i + 2 $. Therefore, there is at least one subsequence $v'h'$ between the $h'v'$ subsequences corresponding to $(x_i, y_i)$ and $(x_{i+1}, y_{i+1})$.
\epr

Figure~\ref{fig:gridWalk} displays an example of a grid walk from $(0,0)$ to $(6,5)$. The endpoints of $h'v'$ subsequences form the strict chain $\{(1,1), (3,3), (4,4), (6,5)\}$. This chain is not maximal, since $(2,2)$ could be added to the chain. This element is the endpoint of a $v'h'$ subsequence in the walk, namely $(1,1),(1,2), (2,2)$. The dashed lines mark the $h'v'$ subsequences in the walk. These subsequences could be replaced by a diagonal move $d'$, playing a similar role to the $d$ in the words of three letters described in Section~\ref{sssec:alignWords}. The walk in Figure~\ref{fig:gridWalk} could be described by the word $d'vhd'd'hd'$.

\begin{remark}\label{rem:walkConvention}
In Proposition~\ref{prop:walkStrictChain}, we associate a strict chain to a walk by focusing on the endpoints of $h'v'$ subsequences. Alternatively, we could focus on the endpoints of $v'h'$ subsequences. They define another strict chain which is maximal iff the walk has no $h'v'$ subsequence disjoint from the $v'h'$ subsequences. In Figure~\ref{fig:gridWalk}, the strict chain associated with endpoints of $v'h'$ subsequences is $\{(2,2), (4,3),(5,4)\}$. There are two $h'v'$ subsequences disjoint from the $v'h'$ subsequences: $\{(0,0), (1,0), (1,1)\}$ and $\{(5,4),(6,4),(6,5)\}$. Two elements can be added to the strict chain $\{(2,2), (3,3),(4,4)\}$, namely, $(1,1)$ and $(6,5)$.  Note also that any strict chain can be associated with a walk by adopting one or the other convention. For instance, the strict chain represented by the blue bullet points in Figure~\ref{fig:gridWalk} can alternatively be obtained by applying the other convention to another walk: the walk obtained by replacing each $h'v'$ subsequence by a $v'h'$ subsequence and conversely for the subsequences $v'h'$ disjoint from the $h'v'$ subsequences. The walk $(0,0), (0,1), (1,1), (2,1), (2,2), (2,3), (3,3),$ $(3,4), (4,4), (4,5), (5,5), (6,5)$ determines the strict chain represented by the blue bullet points in Figure~\ref{fig:gridWalk} if we consider the endpoints of the $v'h'$ subsequences of the walk.
\end{remark}

\begin{remark}\label{rem:hPrimeVPrime}
The moves represented by the letters $h$ (\resp $v$) and $h'$ (\resp $v'$) in Figures~\ref{fig:strChainWord} and~\ref{fig:gridWalk} are the same. The points in the grids are just numbered in different ways in these figures. Therefore, the grid walks in Figure~\ref{fig:strChainWord} start from the NW corner and end up in the SE corner while in Figure~\ref{fig:gridWalk}, they start from the SW corner and end up in the NE corner.
\end{remark}

\subsubsection{Walks and antichains}\label{sssec:walkAntichain}
Consider a walk starting from $(0,m_2)$ and ending up in $(m_1,0)$. Moves $h'$ (\resp $v'$) increase (\resp decrease) the first (\resp second) coordinate by one unit. A walk in \mm{} is described by a sequence of $h'$ and $v'$ moves. The endpoints of $h'v'$ subsequences form an antichain in \mm. Such an antichain is maximal iff the walk has no $v'h'$ subsequence disjoint from the $h'v'$ subsequences. This result is established as Proposition~\ref{prop:walkStrictChain}. A similar result holds if we define the elements of an antichain as the endpoints of $v'h'$ subsequences of a walk.

\begin{figure}[h!!!]
  \centering
  \includegraphics[width=7cm]{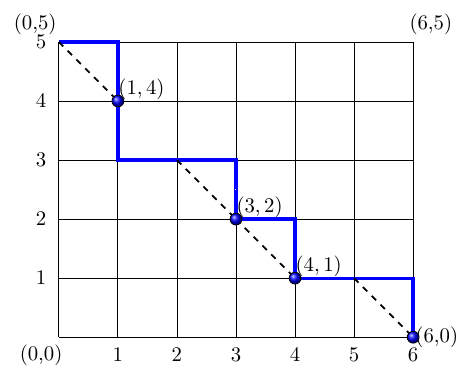}
  \caption{A grid walk from $(0,5)$ to $ (6,0)$. The antichain $\{(1,4), (3,2), (4,1), (6,0)\}$ is not maximal because it has a $v'h'$ subsequence}\label{fig:walkAntichain}
\end{figure}

\ssk Figure~\ref{fig:walkAntichain} represents the grid walk obtained from that in Figure~\ref{fig:gridWalk} by replacing the second coordinate $x_2$ of its points with its complement to $m_2$, \ie $5-x_2$. The elements of the associated antichains are the endpoints of the $h'v'$ subsequences. The antichain is not maximal; it can be augmented with $(1,3)$, the endpoint of the sole $v'h'$ subsequence that is disjoint from the $h'v'$ subsequences.




\subsection{Summary}\label{ssec:summarySec2}
The present section has established correspondences between antichains  (\resp maximal antichains) in the product of two chains and a number of different objects. Table~\ref{ta:correspondences} summarizes these correspondences.

\begin{table}[h!!!]
  \centering
  \begin{tabular}{|l|l|}
    \hline
    Object & Additional property\\
    \hline
    antichain  & maximal \\
    strict chain & maximal  \\
    alignment of two strings & without alternate skips  \\
    $(h,v,d)$-words without $hv$ &  without $vh$ \\
    set of noses of a step matrix & grid line has no $vh$ \\
    endpoints of $h'v'$ subsequences in a walk   & without $v'h'$ subsequence disjoint \\ & \hspace{1cm}  from the $h'v'$ subsequences  \\
    \hline
  \end{tabular}
  \caption{Correspondences}\label{ta:correspondences}
\end{table}

We established such correspondences with most mathematical objects listed in relation with Sequence A171155 in \citet{OEIS2024} up to September 1, 2024.  In order to prove these correspondences, we introduced some new (as far as we know) mathematical objects such as the \emph{augmentation matrix} of an antichain or of a strict chain (Sections~\ref{sssec:augmentMatrixAntich} and~\ref{sssec:StrChainMatrices}. We reused the notion of \emph{nose} of a step matrix (which plays an important role in the numerical representation theory of \emph{semiorders} \citep{Pirlot1990,BalofDoignonFiorini13,BouyssouPirlot2021semiordersCountable,BouyssouPirlot2021semiordersGeneral}, associating them with elements of a strict chain (Proposition~\ref{prop:strChainSWmatrix}), with $d$-letters of words representing certain alignments, or with diagonal $d$-moves in grid lines (Section~\ref{sssec:antichWords}).

\ssk There is a mathematical object mentioned in relation with OEIS sequence A171155 for which we have no correspondence, namely, ``the number of walls of height 1 in bargraphs with semiperimeter $n \geq 2$'' \citep[see][for definitions]{BlecherBrennanKnopfmacher2017}. Bargraphs are a special type of lattice paths. A single bargraph may have several walls of height one (or none). Hence, if a  correspondence does exist with the sort of lattice paths or the sort of words we considered, it cannot be one-to-one. We do not investigate this issue further here.

%
%


\newpage

\section{Counting maximal antichains in products of two chains}\label{sec:counting}
The number of antichains is well-known. We recall one way of enumerating them in Section~\ref{ssec:numbAntich}. Counting maximal antichains is more difficult. In this section, we provide different ways of counting antichains (\resp maximal antichains) in products of two chains, by using some of the correspondences described in Section~\ref{sec:correspondences}.

\subsection{Number of antichains}\label{ssec:numbAntich}

We denote the number of  antichains in $[m_{1}] \times [m_{2}]$ by $d_E(m_1, m_2)$. Clearly, $d_E(m_1,m_2) = d_E(m_2,m_1)$. This number is equal to the number of walks in $\mathbb{Z}^{2}$ from $(0,m_{1})$ to $(m_{2},0)$. Any such walk has $m_1$ horizontal steps $h'$ and $m_2$ vertical steps $v'$. Therefore, the number of such walks is
\begin{equation}\label{eq:numbAntich}
  d_E(m_1, m_2) = {{m_1+m_2} \choose {m_1}}
\end{equation}
When $m_1=m_2$, this is sequence A000984 in \cite{OEIS2024}. One can find there a large number of mathematical objects that can be counted by the above formula with $m_1 = m_2 = m$. 

\subsection{Number of maximal antichains}

We denote by $d_F(m_1, m_2)$ the number of maximal antichains in \mm. This is sequence A180091 in \citet{OEIS2024}. As for antichains, this function is symmetric in its two variables. 

\subsubsection{Heinz's expression}

An expression (stated by Alo\"is P. Heinz) for the number $d_F(m,m)$ can be found in the On-line Encyclopedia of Integer Sequences \cite{OEIS2024}, Sequence A171155:
for all $m \in \NN$,
\begin{align}\label{eq:Heinz}
d_{F}(m,m) =  \frac{1}{m}[&(4m-3)d_{F}(m-1,m-1)-(2m-5)d_{F}(m-2,m-2)\nonumber \\
 & +d_{F}(m-3,m-3)-(m-3)d_{F}(m-4,m-4)].
\end{align}
This expression is computationally efficient, but, unfortunately, no proof is available\footnote{\citet{Heinz_private_com} wrote us
that he most probably found this formula by using computer programs. He had no explicit proof available.}.
%

\subsubsection{An explicit expression}
The following formula is obtained by counting the $(h,v,d)$-words of $2m$ letters without $hv$ nor $vh$ subsequences. By Proposition~\ref{prop:antichWords}, this number is equal to the number of maximal antichains.

\msk
\begin{proposition}\label{prop:explicitFormula} For $m \in \mathbb{N}$,
\begin{align}\label{eq:explicitExpr}
  d_F(m,m) = & 1 + 2 \sum_{k=1}^{\lfloor\frac{2m-1}{3}\rfloor} \sum_{t=1}^{k}{m-k - 1 \choose {\lfloor\frac{t}{2}\rfloor}} {m-k - 1 \choose { \lfloor\frac{t -1}{2}\rfloor}} {2m-k-t \choose k-t} \nonumber\\
  & + 2 \sum_{k=\lfloor\frac{2m-1}{3}\rfloor +1}^{m-1} \sum_{t=1}^{2m-2k-1}{m-k - 1 \choose {\lfloor\frac{t}{2}\rfloor}} {m-k - 1 \choose { \lfloor\frac{t -1}{2}\rfloor}} {2m-k-t \choose k-t}.
\end{align}
\end{proposition}

\bpr
See Appendix \ref{app:sec:explicitExpression} for a proof of this formula.
\epr

We illustrate the use of Formula~\eqref{eq:explicitExpr} for $m=7$ in Appendix~\ref{app:ssec:cases678}. This expression is probably not very efficient from a computational viewpoint. Note that a similar expression has been provided by Thomas Baruchel (November 9, 2014) on \citet{OEIS2024}, seq. A171155.

\subsubsection{A double recurrence for $d_F(m_1, m_2)$}\label{sssec:doubleRec}

Let $d^{h}_{F}(m_{1},m_{2})$ denote the number of walks in $\mathbb{Z}^{2}$ from $(0,m_{1})$ to $(m_{2},0)$ such that the first step is $h$. Such walks consist of $h, v$ and $d$ moves with the constraint that an $h$ move cannot be followed by a $v$ move and a $v$ move cannot be followed by a $h$ move.

\msk
\begin{proposition}\label{prop:doubleRecur} For $m_1, m_2 \in \mathbb{N}$,
\begin{equation}\label{eq:doubleRecur}
\left\{  \begin{array}{l}
    d_{F}(m_{1},m_{2}) =  d^h_{F}(m_{1},m_{2}) + d^h_{F}(m_{2},m_{1}) + d_{F}(m_{1}-1,m_{2}-1); \\
     d^{h}_{F}(m_{1},m_{2}) =  d_{F}(m_{1}-1,m_{2}) - d_{F}(m_{1}-1,m_{2}-1) + d^{h}_{F}(m_{1}-1,m_{2}-1).
  \end{array}
\right.
\end{equation}
\end{proposition}

\bpr
Proposition~\ref{prop:doubleRecur} is established in Appendix~\ref{app:sec:proofDoubleRecur}.
\epr

This double recurrence is computationally less (\resp more) efficient than Heinz's expression (\resp Expression \eqref{eq:explicitExpr}). In Appendix~\ref{app:sec:valuesDoubleRecur}, Table \ref{app:ta:doubleRecur8df} (\resp Table~\ref{app:ta:doubleRecur8dfh}) displays $d_{F}(m_{1},m_{2})$ (\resp $d^h_{F}(m_{1},m_{2})$) up to $m_1, m_2 = 8$.

\subsubsection{A simple recurrence for $d_F(m_1, m_2)$}\label{sssec:simpleRecur}

The following recurrence relation holds for $d_F(m_1, m_2)$. It is due to \citet{Covington2004} who gives a succinct proof of the result. We provide a more detailed proof below.
\begin{proposition}\label{prop:simpleRecurr}
For all $m_{1},m_{2} \in \NN$,
\begin{equation}\label{eq:simpleRecur}
d_{F}(m_{1},m_{2}) = d_{F}(m_{1}-1,m_{2}-1) + \sum_{i=0}^{m_{1}-2} d_{F}(i,m_{2}-1) + \sum_{i=0}^{m_{2}-2} d_{F}(m_{1}-1,i).
\end{equation}
\end{proposition}
\bpr
Let $A$ be the set $\{(1,i) : i \in [m_{2}] \} \cup \{ (i,m_{2}) : i \in [m_{1}] \}$
Let us first prove that
a maximal antichain in $[m_{1}] \times [m_{2}]$ necessarily intersects $A$. Suppose for contradiction that $s$ is a maximal antichain not intersecting $A$. Notice that $(1,m_{2}) \in A$ is incomparable to all elements in $s$.
Then the set $s \cup \{(1,m_{2})\}$ is also an antichain and, hence, $s$ is not maximal. This proves our claim. Since all elements of $A$ are comparable, each maximal antichain in $[m_{1}] \times [m_{2}]$ contains only one element of $A$.
Moreover, it is easy to see that each element of $A$ is contained in at least one maximal antichain.

\ssk Since each maximal antichain contains one  and only one element of $A$, we can partition the set of maximal antichains according to the element of $A$ that they contain. This partition has $m_{1}+m_{2}-1$ elements.
 In order to count the number of maximal antichains, we can add the number of maximal antichains in each element of the partition. Consider any $a \in A$.
 There are three exclusive possible cases.
 \begin{itemize}

 \item $a= (1, m_{2})$. The set of elements incomparable to $(1,m_{2})$ is $[2,m_{1}] \times [1,m_{2}-1]$, where $[a,b] = \{i \in \NN : a \leq i \leq b\}$. The number of maximal antichains containing $a$ is therefore the number of maximal antichains in the poset $[2,m_{1}] \times [1,m_{2}-1]$ which is also the number of maximal antichains in the poset $[m_{1}-1] \times [m_{2}-1]$, that is $d_{F}(m_{1}-1,m_{2}-1)$.  This corresponds to the first term in \eqref{eq:simpleRecur}.

  \item $a= (a_{1},m_{2})$ with $a_{1} \in [2,m_{1}]$. The set of elements incomparable to $(a_{1},m_{2})$ is $[a_{1}+1,m_{1}] \times [1,m_{2}-1]$. The number of maximal antichains containing $a$ is therefore the number of maximal antichains in the poset $[a_{1}+1,m_{1}] \times [1,m_{2}-1]$ which is also the number of maximal antichains in the poset $[m_{1}-a_{1}] \times [m_{2}-1]$. Depending on $a_{1}$, this  varies from $d_{F}(0,m_{2}-1)$ to $d_{F}(m_{1}-2,m_{2}-1)$. This corresponds to the second term in \eqref{eq:simpleRecur}.

\item $a= (1,a_{2})$ with $a_{2} \in [m_{2}-1]$. The set of elements incomparable to $(1,a_{2})$ is $[2,m_{1}] \times [1,a_{2}-1]$. The number of maximal antichains containing $a$ is therefore the number of maximal antichains in the poset $[2,m_{1}] \times [1,a_{2}-1]$ which is also the number of maximal antichains in the poset $[m_{1}-1] \times [a_{2}-1]$. Depending on $a_{2}$, this  varies from $d_{F}(m_{1}-1,0)$ to $d_{F}(m_{1}-1,m_{2}-2)$. This corresponds to the third term in \eqref{eq:simpleRecur}. \epr

 \end{itemize}

\section{Asymptotic results}

This section uses Heinz's result (\ie Eq.\ \eqref{eq:Heinz}) for which no proof is available.

\medskip\noindent

The ratio $d_{F}(m+1,m+1)/d_{F}(m,m)$ lies between 3 and 3.2 for $m \in \{5,6,7\}$ (see Table~\ref{app:ta:doubleRecur8df}). More generally, for all $l \geq 3$ and $m \in \{l-2, l-1,l\}$, the ratio $d_{F}(m+1,2)/d_{F}(m,2)$ lies between some reals $a_{l}$ and $b_{l}$.
Put differently, for all $l \geq 3$ and $m \in \{l-2, l-1,l\}$,
$$d_{F}(m+1,m+1)/b_{l} \leq d_{F}(m,2) \leq d_{F}(m+1,2)/a_{l}.$$
 So, using Heinz's expression, we can write an upper bound for $md_{F}(m,2)$, for all $m \geq 3$:
\begin{align*}
md_{F}(m,2) & = &  (4m-3)d_{F}(m-1,2)-(2m-5)d_{F}(m-2,2) \\
 & & +d_{F}(m-3,2)-(m-3)d_{F}(m-4,2) \\
  & \leq & (4m-3)d_{F}(m-1,2) - (2m-5)d_{F}(m-1,2)/a_{l} \\
 & & +d_{F}(m-1,2)/b^{2}_{l} - (m-3)d_{F}(m-1,2)/a^{3}_{l}  \\
  & = & \left( (4m-3) - \frac{2m-5}{a_{l}} + \frac{1}{b^{2}_{l} } - \frac{m-3}{a^{3}_{l}} \right)  \  d_{F}(m-1,2).
\end{align*}
This implies
$$ \frac{d_{F}(m,2)}{d_{F}(m-1,2)} \leq \frac{4m-3}{m} - \frac{2m-5}{m a_{l}  } + \frac{1}{m b^{2}_{l}} - \frac{m-3}{m a^{3}_{l}}$$
and
$$\lim_{m \to \infty} \frac{d_{F}(m,2)}{d_{F}(m-1,2)} \leq 4 - \frac{2}{ a_{l}} - \frac{1}{a^{3}_{l}}.$$
This expression provides the tightest upper bound for $d_{F}(m,2)/d_{F}(m-1,2)$ when $a_{l}$ is as small as possible. This occurs when
$$a_{l} = 4 - \frac{2}{ a_{l}} - \frac{1}{a^{3}_{l}}.$$
This equation has two real roots: 1 and
\begin{equation*}
\label{root}
\rho =\frac{1}{3} \left( 3 + \sqrt[3]{54-6\sqrt{33}} + \sqrt[3]{6(9+\sqrt{33})} \right) .
\end{equation*}
The first root  is not a solution to our problem because it would imply that $d_{F}(m,2)$ does not grow with $m$. So, the only feasible solution is $\rho$.

If we write a lower bound for $md_{F}(m,2)$ and follow the same reasoning as above, we obtain also $\rho$ for $b_{l}$. We can therefore conclude that
$$   \lim_{m \to \infty} \frac{d_{F}(m,2)}{d_{F}(m-1,2)} = \rho \approx 3.38.$$

For the number of antichains, we can also compute the ratio $d_{E}(m+1,2)/d_{E}(m,2)$. It is equal to
$$\frac{ {2(m+1) \choose m+1} }{{2m \choose m}} = \frac{(2m+2) (2m+1)}{(m+1)^{2}}.$$
In the limit, this ratio is equal to 4.
We can therefore conclude that
$$\lim_{m \to \infty} \frac{d_{F}(m,m)}{d_{E}(m,m)} = 0.$$
Provided Heinz's formula is correct, this shows that the number of maximal antichains is much less than the number of antichains, for large $m$.

\begin{appendix}
\section*{Appendix}

\section{Proof of Proposition~\ref{prop:explicitFormula}}\label{app:sec:explicitExpression}

Proposition~\ref{prop:antichWords} implies that the number of maximal antichains $d_F(m,m)$ in a product of two chains $[m]\times [m]$ is equal to the number of words involving three letters $d,h,v$ with the constraint that $hv$ and $vh$ are forbidden sequences of two consecutive letters. Such words involve  $k$ letters $d$ and an equal number $l$ of $h$ and $v$ letters with $2k + 2l =2m$ or $k+l=m$. How many words of that type are there for fixed $m$?

\ssk We distinguish (and will count) the words by grouping them in subsets involving the same numbers of each of the three letters. A subset of this type is described by a formula such as $k \times d + l \times h + l \times v$.

\begin{example} Let us consider for instance the case $m=7$. The following subsets of words may be distinguished:
\begin{itemize}
  \item $7 d$
  \item $6 d + 1 h + 1 v$
  \item $5 d + 2 h + 2 v$
  \item $4 d + 3 h + 3 v$
  \item $3 d + 4 h + 4 v$
  \item $2 d + 5 h + 5 v$
  \item $1 d + 6 h + 6 v$
\end{itemize}
There are no legal words involving no $d$ letter. Also it is easy to count the words of the first type $7d$: there is only one. For the six other types, there is a symmetry between words in which the first non $d$ letter is $h$ or is $v$. We shall count those words in which the first non $d$ letter is $h$ and double the count. Therefore, the counts of the 6 latter types are even numbers. Since there is only one word of the first type, the total count is always an odd number.
\end{example}

We use the following idea for counting words of each type. Ignoring the $d$ letters at first glance, we count the possible sequences of $h$ and $v$ letters in a word. We call such a sequence a \textit{schema}. By hypothesis, we count only the schemas starting with $h$ letter and we double the count afterwards. We call \textit{transition} a subsequence of consecutive $hv$ or $vh$. In the complete word, each transition must be separated by a $d$ letter. Therefore, the number of transitions in a word must be at most the number $k$ of $d$ letters. This number is also bounded by the number of transitions that can be obtained using $l$ letters $h$ and $l$ letters $v$, i.e., $2 l -1$. The number of transitions $t$ in a word thus satisfies
\begin{equation}\label{eq:numberTransitions}
  t \leq \min(k, 2l-1).
\end{equation}
In the middle of each transition, we need to put a letter $d$. We then remain with $k-t$ letters $d$ and we have to count the manners they can be inserted in the word. This is the counting strategy we apply. We illustrate it with  two examples.

\begin{example}
Let us take for example the type $4d + 3h + 3v$ for $m=7$. There may be up to 4 transitions since $4 = \min(k, 2l-1) = \min(4, 5)$. We count the schemas and the ways of inserting $d$'s in the schemas for $t=1, 2, 3, 4$ transitions. Remember that we count the schemas starting with $h$ and then double the count.
\begin{description}
  \item[1 transition] The 3 $h$ must be before the 3 $v$; there is only one schema: $hhhvvv$. One $d$ is used to separate the only transition: $hhhdvvv$. Then there are 7 positions in which the remaining 3 $d$'s can be inserted. Several ones can be inserted in the same position. The formula to be used is that of combinations with repetitions. We have to chose 3 chocolates from 7 sorts. The number is $ {7+3-1 \choose 3} = {9 \choose 3} = 84$.
  \item[2 transitions] To count the schemas, we have to count the number of ways of inserting 1 $h$ and 2 $v$'s into $hvh$ without creating additional transitions. Thus the $h$ must be inserted together with the $h$'s and the $v$'s with the $v$. There are 2 possibilities for inserting $h$ and 1 for the $v$'s. In each of the two possible schemas $hhvvvh$ and $hvvvhh$, we need 2 $d$'s for separating the two transitions. The number of ways of inserting the two remaining $d$'s is $ {7+2-1 \choose 2} = {8 \choose 2} = 28$. This number multiplied by the number of schemas makes 56.
  \item[3 transitions] To insert one $h$ and one $v$ into $hvhv$, we have $2 \times 2 = 4$ possibilities hence 4 schemas. We use three $d$'s for separating the transitions. One $d$ remains to be inserted; there are $ {7+1-1 \choose 1} = {7 \choose 1} =7$ ways of doing this. The count for this case is thus $4 \times 7 = 28$.
  \item[4 transitions] We have to insert one $v$ into $hvhvh$; there are thus 2 schemas. All four $d$'s are used to separate the transitions. The total count is thus 2.
\end{description}
Summing up, we have $(84 + 56 + 28 + 2) \times 2 = 340$ words satisfying the formula $4d + 3h + 3v$.
\end{example}

\begin{example} Let us take another example: $5d + 2h + 2v$. The number of transitions here is bounded by $\min(5, 4-1) = 3$.
\begin{description}
  \item[1 transition] There is only one way of inserting an additional $h$ and an additional $v$ in the sequence $hv$. Thus, only one schema: $hhvv$. One $d$ is used; we need to count the ways of inserting the remaining four $d$'s (choosing 4 chocolates). There are 5 possible positions (sorts of chocolates). The number of possibilities is ${5+4-1 \choose 4} = {8 \choose 4}= 70$.
  \item[2 transitions] There is one way of inserting an additional $v$ into $hvh$. So, one schema. The number of possibilities for inserting the three remaining $d$'s is ${5 + 3 -1 \choose 3} = 35$.
  \item[3 transitions] There is only one schema: $hvhv$. It uses 3 out of the 5 $d$'s. There are ${5+2-1 \choose 2} = 15$ ways of inserting the remaining 2 $d$'s.
\end{description}
The total count of the words of the type $5d + 2h +2v$ is $(70 + 35 + 15) \times 2 = 240$.
\end{example}

\ssk We now apply our counting strategy in the general case.  We want to count the number of words of type $kd + lh + lv$. The number of transitions is bounded by $\min(k, 2l-1)$.

\ssk We distinguish the cases in which $t$ is even or odd.

\begin{description}
  \item[$t = 2t'$ is even] The alternate sequence of $h$'s and $v$'s that determine the $2t'$ transitions involves $t'+1$ letters $h$ and $t'$ letters $v$. We have to count the number of ways of inserting the remaining $l-t'-1$ letters $h$ and the remaining $l-t'$ letters $v$. For the $h$'s, we have to insert (possibly with repetitions) $l-t'-1$ letters in the $t'+1$ positions already held by the $h$'s. There are ${(l-t'-1) + (t'+ 1) - 1 \choose (l-t'-1)} = {l - 1 \choose t'}$ ways of doing that (note that the chocolates here are the $l-t'-1$ remaining $h$'s, to be chosen among $t'+1$ sorts of chocolates). Similarly, there are ${(l-t') + (t') - 1 \choose (l-t')} = {l - 1 \choose t' -1}$ ways of inserting the remaining $t'$ letters $v$. The number of schemas is thus
      $$
      {l - 1 \choose t'}\times {l - 1 \choose t' -1}.
      $$

  \item[$t = 2t'-1 $ is odd] The alternate sequence of $h$'s and $v$'s that determine the $2t'-1$ transitions involves $t'$ letters $h$ and $t'$ letters $v$. The number of schemas in this case is thus:
      $$
      {l - 1 \choose t'-1}\times {l - 1 \choose t' -1}.
      $$
\end{description}
The following expression for the number of schemas covers both the cases of $t$ even and odd:
      $$
      {l - 1 \choose {\lfloor\frac{t}{2}\rfloor}}\times {l - 1 \choose { \lfloor\frac{t -1}{2}\rfloor}}.
      $$

For $t$ transitions, we use $t$ out of the $k$ letters $d$ that are available. The number of ways of inserting the remaining $k-t$ letters $d$ into the $2l+1$ positions in the schemas is ${(2l+1) + (k-t) -1 \choose k-t} = {2l+k-t \choose k-t}$.

The number of words of type $kd+lh+lv$, with $t$ transitions ($1 \leq t \leq \min(k, 2l-1) $) is

      $$
      2 \times {l - 1 \choose {\lfloor\frac{t}{2}\rfloor}}\times {l - 1 \choose { \lfloor\frac{t -1}{2}\rfloor}} \times {2l+k-t \choose k-t}.
      $$
Since $l, k$ and $m$ are linked by $k+l =m$, we eliminate $l$. The above formula becomes:
      $$
      2 \times {m-k - 1 \choose {\lfloor\frac{t}{2}\rfloor}}\times {m-k - 1 \choose { \lfloor\frac{t -1}{2}\rfloor}} \times {2m-k-t \choose k-t},
      $$
for $1 \leq t \leq \min(k, 2m-2k-1)$.

For $k \leq 2l-1 = 2m-2k-1$, the number of transitions $t$ varies from 1 to $k$. For larger values of $k$, $t$ is bounded by $2l-1$. Solving the inequality $k \leq 2m-2k-1$ yields $k \leq \frac{2m-1}{3}$. So, the formula for $d_F(m,m)$ reads as follows:
\begin{align*}\label{eq:DFm2}
  d_F(m,m) = & 1 + 2 \sum_{k=1}^{\lfloor\frac{2m-1}{3}\rfloor} \sum_{t=1}^{k}{m-k - 1 \choose {\lfloor\frac{t}{2}\rfloor}} {m-k - 1 \choose { \lfloor\frac{t -1}{2}\rfloor}} {2m-k-t \choose k-t} \nonumber\\
  & + 2 \sum_{k=\lfloor\frac{2m-1}{3}\rfloor +1}^{m-1} \sum_{t=1}^{2m-2k-1}{m-k - 1 \choose {\lfloor\frac{t}{2}\rfloor}} {m-k - 1 \choose { \lfloor\frac{t -1}{2}\rfloor}} {2m-k-t \choose k-t}
\end{align*}




\clearpage
\section{Using Formula~\eqref{eq:explicitExpr} for $m=7$}\label{app:ssec:cases678}
We use Formula~\eqref{eq:explicitExpr} to compute 
$d_F(7,7) = 817$ (Table \ref{ta:DF72})
.



\begin{table}[h!!!]
  \centering
  \begin{tabular}{|l|l|l|r|r|}
    \hline
    $k$ & $t$ & \# words & count & total $\times 2$  \\
    \hline
    1 & 1 & ${5 \choose 0} {5 \choose 0}{12 \choose 0}$ &  1 & 2 \\
    \hline
    2 & 1 & ${4 \choose 0} {4 \choose 0}{11 \choose 1}$ & 11 &  \\
      & 2 & ${4 \choose 1} {4 \choose 0}{10 \choose 0}$ &  4 &  30 \\
    \hline
    3 & 1 & ${3 \choose 0} {3 \choose 0}{10 \choose 2}$ & 45 &    \\
      & 2 & ${3 \choose 1} {3 \choose 0}{ 9 \choose 1}$ & 27 &    \\
      & 3 & ${3 \choose 1} {3 \choose 1}{ 8 \choose 0}$ &  9 & 162\\
    \hline
    4 & 1 & ${2 \choose 0} {2 \choose 0}{ 9 \choose 3}$ & 84 & \\
      & 2 & ${2 \choose 1} {2 \choose 0}{ 8 \choose 2}$ & 56 & \\
      & 3 & ${2 \choose 1} {2 \choose 1}{ 7 \choose 1}$ & 28 & \\
      & 4 & ${2 \choose 2} {2 \choose 1}{ 6 \choose 0}$ &  2 & 340 \\
    \hline
    5 & 1 & ${1 \choose 0} {1 \choose 0}{ 8 \choose 4}$ & 70 & \\
      & 2 & ${1 \choose 1} {1 \choose 0}{ 7 \choose 3}$ & 35 & \\
      & 3 & ${1 \choose 1} {1 \choose 1}{ 6 \choose 2}$ & 15 & 240 \\
    \hline
    6 & 1 & ${0 \choose 0} {0 \choose 0}{ 7 \choose 5}$ & 21 & 42\\
    \hline
    7 &   &                                             &    &  1\\
    \hline \hline
      &   &                                             &    & 817\\
    \hline
  \end{tabular}
  \caption{Computation of $d_F(7,7)$ using formula~\eqref{eq:explicitExpr}. The number of transitions is limited by $2m-2k-1 = 14 - 2k -1$ as soon as $k > \lfloor \frac{2m-1}{3}\rfloor = 4$}\label{ta:DF72}
\end{table}

\section{Proof of Proposition~\ref{prop:doubleRecur}}\label{app:sec:proofDoubleRecur}
For all $m,n \in \NN$, let $\dfh(m,n)$ (\resp $\dfv(m,n)$, $\dfd(m,n)$) denote the number of walks starting with a $h$ (\resp $v$, $d$) move. Note that for all $m, n \in \NN$,   $\df(m, n) = \df(n, m)$ and $\dfh(n,m) = \dfv(n,m)$. We have:
$$
\df(\un,\deu) = \dfh(\un, \deu) + \dfv(\un, \deu) + \dfd(\un, \deu).
$$
Since $\dfv(\un, \deu) = \dfh(\deu, \un)$ and $\dfd(\un, \deu) = \df(\un-1, \deu-1)$, we get:
\begin{equation}\label{eq:df}
  \df(\un,\deu) = \dfh(\un, \deu) + \dfh(\deu, \un) + \df(\un-1, \deu-1).
\end{equation}
We may also write:
\begin{align}\label{eq:dfh}
  \dfh(\un, \deu)  & = \df(\un-1, \deu-1) - \dfv(\un-1, \deu) \nonumber \\
   & = \df(\un-1, \deu-1) - \dfh(\deu, \un-1).
\end{align}
Applying the same idea to $\dfh(\deu, \un-1 $, we get:
\begin{align}\label{eq:dfhDeuUn}
  \dfh(\deu, \un -1)  & = \df(\deu-1, \un-1) - \dfv(\deu-1, \un-1) \nonumber \\
   & = \df(\deu-1, \un-1) - \dfh(\un-1, \deu-1).
\end{align}
Substituting \eqref{eq:dfhDeuUn} into \eqref{eq:dfh} yields
\begin{equation}\label{eq:formulaDfh}
    \dfh(\un, \deu) =   \df(\un-1, \deu-1) -   \df(\deu -1, \un -1) +   \dfh(\un -1, \deu -1).
\end{equation}

This establishes Proposition~\ref{prop:doubleRecur}.

\section{Values given by the double recurrence~\eqref{eq:doubleRecur}}\label{app:sec:valuesDoubleRecur}
Using the double recurrence~\eqref{eq:doubleRecur}, we computed the values of $\df(\un,\deu)$ and $\dfh(\un, \deu)$ for $\un, \deu$ up to 8 displayed in Tables \ref{app:ta:doubleRecur8df} and \ref{app:ta:doubleRecur8dfh}.

\begin{table}[h!!!]
  \centering
  \begin{tabular}{|l|r|r|r|r|r|r|r|r|}
    \hline
    $\df$ & $\deu=1$ & $\deu=2$ & $\deu=3$ & $\deu=4$ & $\deu=5$ & $\deu=6$ & $\deu=7$ & $\deu=8$ \\ \hline
    $\un=1$ & 1 & 2 & 3 & 4 & 5 & 6 & 7 & 8 \\
    $\un=2$ & 2 & 3 & 5 & 8 & 12 & 17 & 23 & 30 \\
    $\un=3$ & 3 & 5 & 9 & 15 & 24 & 37 & 55 & 79 \\
    $\un=4$ & 4 & 8 & 15 & 27 & 46 & 75 & 118 & 180 \\
    $\un=5$ & 5 & 12 & 24 & 46 & 83 & 143 & 237 & 380 \\
    $\un=6$ & 6 & 17 & 37 & 75 & 143 & 259 & 450 & 755 \\
    $\un=7$ & 7 & 23 & 55 & 118 & 237 & 450 & 817 & 1429 \\
    $\un=8$ & 8 & 30 & 79 & 180 & 380 & 755 & 1429 & 2599 \\
    \hline
  \end{tabular}
  \caption{Values of $\df(\un,\deu)$ for $\un, \deu = 1, \ldots, 8 $}\label{app:ta:doubleRecur8df}
\end{table}

\begin{table}[h!!!]
  \centering
  \begin{tabular}{|l|r|r|r|r|r|r|r|r|}
    \hline
    $\dfh$ & $\deu=1$ & $\deu=2$ & $\deu=3$ & $\deu=4$ & $\deu=5$ & $\deu=6$ & $\deu=7$ & $\deu=8$ \\ \hline
    $\un=1$ & 0 & 0 & 0 & 0 & 0 & 0 & 0 & 0 \\
    $\un=2$ & 1 & 1 & 1 & 1 & 1 & 1 & 1 & 1 \\
    $\un=3$ & 2 & 2 & 3 & 4 & 5 & 6 & 7 & 8 \\
    $\un=4$ & 3 & 4 & 6 & 9 & 13 & 18 & 24 & 31 \\
    $\un=5$ & 4 & 7 & 11 & 18 & 28 & 42 & 61 & 86 \\
    $\un=6$ & 5 & 11 & 19 & 33 & 55 & 88 & 136 & 204 \\
    $\un=7$ & 6 & 16 & 31 & 57 & 101 & 171 & 279 & 441 \\
    $\un=8$ & 7 & 22 & 48 & 94 & 176 & 314 & 538 & 891 \\
    \hline
  \end{tabular}
  \caption{Values of $\dfh(\un,\deu)$ for $\un, \deu = 1, \ldots, 8 $}\label{app:ta:doubleRecur8dfh}
\end{table}
\end{appendix}

\clearpage

\bibliographystyle{plainnat}
\bibliography{AboutMaximalAntichains,main-minimal}
%

\end{document}